\documentclass{amsart}
\usepackage{amssymb}

\theoremstyle{plain}
\newtheorem{theorem}{Theorem}
\newtheorem{lemma}{Lemma}
\theoremstyle{definition}
\newtheorem{problem}{Problem}
\newcommand{\new}[1]{M_3\langle#1\rangle}

\newcommand{\set}[1]{\{#1\}}
\newcommand{\setm}[2]{\{\,#1\mid#2\,\}}
\newcommand{\es}{\varnothing}
\newcommand{\ci}{\subseteq}
\newcommand{\ii}{\cap}
\newcommand{\jj}{\vee}
\newcommand{\mm}{\wedge}
\newcommand{\ol}[1]{\overline{#1}}
\newcommand{\F}[1]{\mathfrak{#1}}
\newcommand{\tbf}{\textbf}
\newcommand{\mbf}{\mathbf}
\newcommand{\gF}{\Phi}
\newcommand{\gQ}{\Theta}
\def\con#1=#2(#3){#1\equiv#2\pod{#3}}
\def\vv<#1>{\langle#1\rangle}

\begin{document}

\title{Proper congruence-preserving extensions of lattices}
\author{G.~Gr\"atzer}
\thanks{The research of the first author was supported by the
        NSERC of Canada.}
\address{Department of Mathematics\\
	  University of Manitoba\\
	  Winnipeg MN, R3T 2N2\\
	  Canada}
\email{gratzer@cc.umanitoba.ca}
\urladdr{http://server.maths.umanitoba.ca/homepages/gratzer.html/}

\author{F.~Wehrung}
\address{C.N.R.S., E.S.A. 6081\\
   D\'epartment de Math\'ematiques\\
	  Universit\'e de Caen\\
	  14032 Caen Cedex\\
   France}
\email{gremlin@math.unicaen.fr}
\urladdr{http://www.math.unicaen.fr/\~{}wehrung}

 \date{February 20, 1998}
\begin{abstract}
We prove that every lattice with more than
one element has a proper congruence-preserving extension.
 \end{abstract}

 \maketitle

\section{Introduction}
Let $L$ be a lattice.  A lattice $K$ is
a \emph{congruence-preserving extension} of $L$, if $K$ is an
extension and every congruence of $L$ has exactly one extension
to $K$.  (Of course, then, the congruence lattice of
$L$ is isomorphic to the congruence lattice of~$K$.)

In \cite{GS95}, the first author and E. T. Schmidt raised the
following question:

Is it true that every lattice $L$ with more than one element
has a proper congru\-ence-preserving extension $K$?

Here \emph{proper} means that $K$ properly contains $L$, that
is, $K - L \neq \es$.

The first author and E. T. Schmidt pointed out in \cite{GS95}
that in the finite case this is obviously true, and they
proved the following general result:

 \begin{theorem}\label{T:old}
 Let $L$ be a lattice.  If there exist a distributive interval
with more than one element in $L$, then $L$ has a proper
congruence-preserving extension $K$.
  \end{theorem}

Generalizing this result, in this paper, we provide a positive
answer to the above question:

 \begin{theorem}\label{T:new}
 Every lattice $L$ with more than one element has a proper
congruence-preserving extension $K$.
  \end{theorem}

\section{Background} Let $K$ and $L$ be lattices.  If $L$ is a
sublattice of $K$, then we call
$K$ an \emph{extension} of $L$.  If $K$ is an extension of $L$ and
$\gQ$ is a congruence relation of $K$, then $\gQ_L$, the
restriction of
$\gQ$ to $L$ is a congruence of $L$.  If the map $\gQ \mapsto
\gQ_L$ is a bijection between the congruences of $L$ and the
congruences of $K$, then we call $K$ a
\emph{congruence-preserving extension} of~$L$.  Observe that
if $K$ a congruence-preserving extension of~$L$, then the
congruence lattice of $L$ is isomorphic to the  congruence
lattice of $K$ in a natural way.

The proof of Theorem~\ref{T:old} is based on the following
construction of E. T. Schmidt \cite{tS62}, summarized below as
Theorem~\ref{T:Schmidt}. (A number of papers utilize this
construction; see, for instance, E. T. Schmidt \cite{tS68},
\cite{tS74} and the recent paper G. Gr\"atzer and E.~T.
Schmidt \cite{GSnew}.)

 Let $L$ be a bounded distributive lattice with bounds $0$ and
$1$, and let $M_3 = \set{o, a, b, c, i}$ be the five-element
nondistributive modular lattice.  Let $M_3[L]$ denote the
poset of triples $\vv<x,y,z> \in L^3$ satisfying the condition
 \begin{equation}
   x \mm y = y \mm z = z \mm x \tag{S}.
 \end{equation}

\begin{theorem}\label{T:Schmidt}\hfill

Let $D$ be a bounded distributive lattice with bounds $0$ and
$1$.
 \begin{enumerate}
 \item $M_3[D]$ is a modular lattice.
 \item The subset
 \[
    \ol{M}_3 = \set{\vv<0, 0, 0>, \vv<1, 0, 0>, \vv<0, 1, 0>,
    \vv<0, 0, 1>, \vv<1, 1, 1>}
 \]
of $M_3[D]$ is a sublattice of $M_3[D]$ and it is
isomorphic to $M_3$.
 \item The subposet $\ol D = \setm{\vv<x, 0, 0>}{x \in D}$ of
$M_3[D]$ is a bounded distributive lattice and it is isomorphic to
$D$; we identify $D$ with
$\ol D$.\label{3}
 \item $\ol{M}_3$ and $D$ generate $M_3[D]$.
  \item Let $\gQ$ be a congruence relation of $D = \ol D$; then
there is a
\emph{unique} congruence $\ol{\gQ}$ of $M_3[D]$ such that
$\ol{\gQ}$ restricted to $\ol D$ is $\gQ$; therefore, $M_3[D]$
is a congruence-preserving extension of~$D$.
 \end{enumerate}
 \end{theorem}

Unfortunately, $M_3[L]$ fails, in general, to produce
a lattice, if $L$ is not distributive.

In this paper, we introduce a variant on the $M_3[L]$
construction, which we shall denote as $\new{L}$.  This lattice
$\new{L}$ is a proper congruence-preserving extension of $L$,
for any lattice $L$ with more than one element, verifying
Theorem~\ref{T:new}.

\section{The construction}
 For a lattice $L$, let us call the triple $\vv<x,y,z> \in L^3$
\emph{Boolean}, if
 \begin{align}
   x  &= (x \jj y) \mm  (x \jj z),\notag\\
   y  &= (y \jj x) \mm  (y \jj z),\tag{B}\\
   z  &= (z \jj x) \mm  (z \jj y)\notag.
 \end{align}

We denote by $\new{L} \ci L^3$ the poset of Boolean triples of
$L$.

Here are some of the basic properties of Boolean triples:

 \begin{lemma}\label{L:basic} Let $L$ be a lattice.\hfill
 \begin{enumerate}
 \item Every Boolean triple of $L$ satisfies \textup{(S)}, so
$\new{L} \ci M_3[L]$.
 \item $\vv<x,y,z> \in L^3$ is Boolean iff there is a triple
$\vv<u,v,w> \in L^3$ satisfying
 \begin{align}
   x  &= u \mm v,\notag\\
   y  &= u \mm w,\tag{R}\\
   z  &= v \mm w\notag.
 \end{align}
 \item\label{I:close} For every triple $\vv<x, y, z> \in L^3$,
there is a smallest Boolean triple $\ol{\vv<x, y, z>} \in L^3$
such that $\vv<x, y, z> \leq \ol{\vv<x, y, z>}$; in fact,
\[
   \ol{\vv<x, y, z>} = \vv<(x \jj y) \mm (x \jj z), (y \jj x)
\mm (y \jj
       z), (z \jj x) \mm (z \jj y) >.
 \]
 \item\label{I:ops}
$\new{L}$ is a lattice with the meet
operation defined as
 \[
   \vv<x_0, y_0, z_0> \mm \vv<x_1, y_1, z_1> =
         \vv<x_0 \mm x_1, y_0 \mm y_1, z_0 \mm z_1>
 \]
and the join operation defined by
 \[
   \vv<x_0, y_0, z_0> \jj \vv<x_1,y_1,z_1> =
       \ol{\vv<x_0 \jj x_1, y_0 \jj y_1, z_0 \jj z_1>}.
 \]
 \item  If $L$ has $0$, then the subposet $\setm{\vv<x, 0, 0>}{x \in
L}$ is a sublattice and it is isomorphic to $L$.

If $L$ has $0$ and $1$, then $\new{L}$ has a spanning
$M_3$, that is, a $\set{0,1}$-sublattice isomorphic to $M_3$,
namely,
 \[
   \set{\vv<0, 0, 0>, \vv<1, 0, 0>, \vv<0, 1, 0>, \vv<0,
       0, 1>, \vv<1, 1, 1>}.
 \]

\item If $\vv<x,y,z>$ is Boolean, then one of the following
holds:
 \begin{enumerate}
 \item the components form a one-element set, so $\vv<x,y,z> =
\vv<a,a,a>$, for some $a \in L$;
 \item the components form a two-element set and
$\vv<x,y,z>$ is of the form $\vv<b,a,a>$, or $\vv<a,b,a>$, or
$\vv<a,a,b>$, for some $a$, $b \in L$, $a < b$.
 \item the components form a three-element set and two
components are comparable and $L$ has two incomparable
elements $a$ and $b$ such that
$\vv<x,y,z>$ is of the form $\vv<a,b,a \mm b>$, or $\vv<a,a \mm
b,b>$, or
$\vv<a \mm b,a,b>$.
 \item the components form a three-element set and the
components are pairwise incomparable and $L$ has an
eight-element Boolean sublattice
$B$ so that the components are the atoms of~$B$.
 \end{enumerate}
 \end{enumerate}
 \end{lemma}

 \begin{proof}\hfill

 (i) If $\vv<x,y,z>$ is Boolean, then
  \begin{align*}
  x \mm y &= ((x \jj y) \mm  (x \jj z)) \mm ((y \jj x) \mm  (y
\jj z))\\
          &= (x \jj y) \mm (y \jj z) \mm  (z \jj x),
 \end{align*}
 which is the upper median of $x$, $y$, and $z$.  So (S) holds.

 (ii) If  $\vv<x,y,z>$ is Boolean, then $u = x \jj y$, $v = x
\jj z$, and $w = y \jj z$ satisfy~(R).  Conversely, if there is
a triple $\vv<u,v,w> \in L^3$ satisfying~(R), then by Lemma
I.5.9 of~\cite{GLT}, the sublattice generated by
$x$, $y$, and $z$ is isomorphic to a quotient of $\F C_2^3$
(where $\F C_2$ is the two element chain) and $x$, $y$, and $z$
are the images of the three atoms of  $\F C_2^3$. Thus
$(x \jj y) \mm  (x \jj z) = x$, the first part of (B).  The other
two parts are proved similarly.

 (iii) For $\vv<x,y,z> \in L^3$, define $u = x \jj y$,
$v = x \jj z$, $w = y \jj z$.  Set $x_1 = u \mm v$,
$y_1 = u \mm w$, $z_1 = v \mm w$.  Then
$\vv<x_1,y_1,z_1>$ is Boolean by (ii) and $\vv<x,y,z> \le
\vv<x_1,y_1,z_1>$ in $L^3$.  Now if $\vv<x,y,z> \le
\vv<x_2,y_2,z_2>$ in
$L^3$ and $\vv<x_2,y_2,z_2>$ is Boolean, then
\begin{align*}
   x_2 &= (x_2 \jj y_2) \mm (x_2 \jj z_2) &&\text{(by (B))}\\
       &\geq  (x \jj y) \mm (x \jj z) &&\text{(by $\vv<x,y,z>
\le
\vv<x_2,y_2,z_2>$)}\\
       &= u \mm v = x_1,
\end{align*} and similarly, $y_2 \geq y_1$, $z_2 \geq z_1$.
Thus $\vv<x_2,y_2,z_2> \geq
\vv<x_1,y_1,z_1>$, and so $\vv<x_1,y_1,z_1>$ is the smallest
Boolean triple containing $\vv<x,y,z>$.

 (iv) $\new{L} \neq \es$; for instance, for all $x \in L$, the
diagonal element $\vv<x, x, x> \in \new{L}$. It is obvious from
(ii) that  $\new{L}$ is meet closed.  By (iii),  $\new{L}$ is
a closure system in $L^3$, from which the formulas of (iv)
follow.

 The proofs of (v) and (vi) are left to the reader.
\end{proof}

 \section{Proof of the theorem} Let $L$ be a lattice with more
than one element.  We identify $x \in L$ with the diagonal
element $\vv<x, x, x> \in \new{L}$, so we regard
$\new{L}$ an extension of $L$.  This is an embedding of $L$ into
$\new{L}$ different from the embedding in
Lemma~\ref{L:basic}.(v).  Moreover, the embedding in
Lemma~\ref{L:basic}.(v) requires that $L$ have a zero, while
the embedding discussed here always works.

 Note that $\new{L}$ is a proper extension; indeed, since $L$
has more than one element, we can choose the elements $a < b$
in $L$.  Then $\vv<a, a, b> \in \new{L}$ but $\vv<a, a, b>$ is
not on the diagonal, so $\vv<a, a, b> \in \new{L} - L$.  In
fact, if $L = \F C_2$, the two-element chain, then this is the
only type of nondiagonal element:
 \[
   \new{\F C_2} = \set{\vv<0,0,0>, \vv<1,0,0>,
     \vv<0,1,0>, \vv<0,0,1>, \vv<1,1,1>}.
 \]

For a congruence $\gQ$ of $L$, let $\gQ^3$ denote the
congruence of
$L^3$ defined component\-wise.  Let $\new{\gQ}$ be the
restriction of
$\gQ^3$ to $\new{L}$.

 \begin{lemma}
 $\new{\gQ}$ is a congruence relation of $\new{L}$.
 \end{lemma}

 \begin{proof}
 $\new{\gQ}$ is obviously an equivalence relation on $\new{L}$.
Since $\new{L}$ is a meet subsemilattice of $L^3$, it is clear
that
$\new{\gQ}$ satisfies the Substitution Property for meets.  To
verify for
$\new{\gQ}$ the Substitution Property for joins, let $\vv<x_0,
y_0, z_0>$, $\vv<x_1, y_1, z_1> \in \new{L}$, let
 \[
   \con \vv<x_0, y_0, z_0>=\vv<x_1, y_1, z_1>(\new{\gQ}),
 \]
 (that is,
\[
 \con x_0=x_1(\gQ),\quad \con y_0=y_1(\gQ), \text{\quad and\quad}
    \con z_0=z_1(\gQ)
 \]
in $L$) and let $\vv<u, v, w> \in \new{L}$.  Set
 \[
   \vv<x_i', y_i', z_i'> = \vv<x_i, y_i, z_i> \jj \vv<u, v, w>
 \]
 (the join formed in $\new{L}$), for $i = 0$, $1$.

Then, using Lemma~\ref{L:basic}.(iii) and (iv) for $x_0 \jj
u$, $y_0
\jj v$, and $z_0 \jj w$, we obtain that
 \begin{align*}
   x_0' &= (x_0 \jj u \jj y_0 \jj v) \mm (x_0 \jj u \jj z_0 \jj
w)\\
        &\equiv (x_1 \jj u \jj y_1 \jj v) \mm (x_1 \jj u \jj
z_1 \jj w) = x_1'\,\,(\new{\gQ}),
 \end{align*}
 and similarly, $\con y_0'=y_1'(\new{\gQ})$,
$\con z_0'=z_1'(\new{\gQ})$, hence
 \[
   \con \vv<x_0, y_0, z_0> \jj \vv<u, v, w>=\vv<x_1, y_1, z_1>
\jj \vv<u, v, w>(\new{\gQ}).
 \]

\vspace{-12pt}
 \end{proof}

Since $L$ was identified with the diagonal of $\new{L}$, it is
obvious that $\new{\gQ}$ restricted to $L$ is $\gQ$.  So to
complete the proof of Theorem~\ref{T:new}, it is sufficient to
verify the following statement:

\begin{lemma}
 Every congruence of $\new{L}$ is of the form $\new{\gQ}$, for a
suitable congruence $\gQ$ of $L$.
 \end{lemma}

 \begin{proof}
 Let $\gF$ be a congruence of $\new{L}$, and let $\gQ$ denote
the congruence of $L$ obtained by restricting $\gF$ to the
diagonal of
$\new{L}$, that is, $\con x=y(\gQ)$ in $L$ iff $\con \vv<x, x,
x>=\vv<y, y, y>(\gF)$ in $\new{L}$. We prove that $\gF =
\new{\gQ}$.

To show that $\gF \ci \new{\gQ}$, let
 \begin{equation}\label{E:1}
   \con \vv<x_0, y_0, z_0>=\vv<x_1, y_1, z_1>(\gF).
 \end{equation} Define
 \begin{align}
   o &= x_0 \mm x_1 \mm y_0 \mm y_1 \mm z_0 \mm
z_1,\label{E:ooo}\\
   i &= x_0 \jj x_1 \jj y_0 \jj y_1 \jj z_0 \jj
z_1\label{E:iii}.
 \end{align}

Meeting the congruence \eqref{E:1} with $\vv<i, o, o>$ yields
 \begin{equation}\label{E:2}
   \con \vv<x_0, o, o>=\vv<x_1, o, o>(\gF).
 \end{equation} Since
 \[
   \vv<x_0, o, o> \jj \vv<o, o, i> = \ol{\vv<x_0, o, i>} =
\vv<x_0, x_0, i>,
 \]
 joining the congruence \eqref{E:2} with $\vv<o, o, i>$ yields
 \begin{equation}\label{E:3}
   \con \vv<x_0, x_0, i>=\vv<x_1, x_1, i>(\gF).
 \end{equation}
 Similarly,
 \begin{equation}\label{E:3new}
   \con \vv<x_0, i, x_0>=\vv<x_1, i, x_1>(\gF).
 \end{equation}
 Now we meet the congruences \eqref{E:3} and
\eqref{E:3new} to obtain
 \begin{equation}\label{E:6}
   \con \vv<x_0, y_0, z_0>=\vv<x_1, y_1, z_1>(\gQ^3)
 \end{equation} in $L^3$, proving that $\gF \ci \new{\gQ}$.

To prove the converse, $\new{\gQ} \ci \gF$, take
 \begin{equation}\label{E:7}
   \con \vv<x_0, y_0, z_0>=\vv<x_1, y_1, z_1>(\new{\gQ})
 \end{equation} in $\new{L}$, that is,
 \begin{alignat*}{2}
    \con x_0&=x_1&&(\gQ),\\
    \con y_0&=y_1&&(\gQ),\\
    \con z_0&=z_1&&(\gQ)
\end{alignat*} in $L$.  Equivalently,
 \begin{alignat}{2}
 \con \vv<x_0, x_0, x_0>&=\vv<x_1, x_1,
x_1>&&(\gF),\label{E:8}\\
 \con \vv<y_0, y_0, y_0>&=\vv<y_1, y_1,
y_1>&&(\gF),\label{E:9}\\
 \con \vv<z_0, z_0, z_0>&=\vv<z_1, z_1, z_1>&&(\gF)\label{E:10}
 \end{alignat} in $\new{L}$.

Now, define $o$, $i$ as in \eqref{E:ooo} and \eqref{E:iii}.
Meeting the congruence
\eqref{E:8} with $\vv<i, o, o>$, we obtain
 \begin{equation}\label{E:11}
   \con \vv<x_0, o, o>=\vv<x_1, o, o>(\gF).
 \end{equation} Similarly, from \eqref{E:9} and \eqref{E:10},
we obtain the congruences
 \begin{alignat}{2}
 \con \vv<o, y_0, o>&=\vv<o, y_1, o>&&(\gF),\label{E:12}\\
 \con \vv<o, o, z_0>&=\vv<o, o, z_1>&&(\gF).\label{E:13}
 \end{alignat} Finally, joining the congruences
\eqref{E:11}--\eqref{E:13}, we get
 \begin{equation}\label{E:14}
   \con \vv<x_0, y_0, z_0>=\vv<x_1, y_1, z_1>(\gF),
 \end{equation} that is, $\new{\gQ} \ci \gF$.  This completes
the proof of this lemma and of Theorem~\ref{T:new}.
 \end{proof}

 \section{Discussion}
 \subsection*{Special extensions}
 We can get a slightly stronger result by requiring that the
extension preserve the zero and the unit, provided they
exist.  To state this result, we need the following concept.

An extension $K$ of a lattice $L$ is \emph{extensive}, provided
that the convex sublattice of $K$ generated by $L$ is $K$.

Note that if $L$ has a zero, $0$, then an extensive extension is
a $\set{0}$-extension (and similarly for the unit, $1$); if $L$
has a zero, $0$, and unit $1$, then an extensive extension is a
$\set{0, 1}$-extension.

 \begin{theorem}\label{T:new-end}
 Every lattice $L$ with more than one element has a proper
congruence-preserving extensive extension $K$.
 \end{theorem}
\begin{proof}
 Indeed, every $\vv<x, y, z> \in \new{L}$ is in the convex
sublattice generated by $L$ since
 \[
   \vv<x \mm y \mm z, x \mm y \mm z, x \mm y \mm z> \leq \vv<x,
y, z>
  \leq \vv<x \jj y \jj z, x \jj y \jj z, x \jj y \jj z>.
 \]

\vspace{-18pt}
 \end{proof}

In Theorem~\ref{T:Schmidt}.(iii), we pointed out that $M_3[D]$
is a congruence-preserving extension of $D = \ol D$, where
$\ol D$ is an ideal of $M_3[D]$.  This raises the question
whether Theorem~\ref{T:new} can be strengthened by requiring
that $L$ be an ideal in $K$.  This is easy to do, if $L$ has a
zero, $0$, since then we can identify $x \in L$ with
$\vv<x, 0, 0> \in \new{L}$.

 \begin{theorem}\label{T:new-ideal}
 Every lattice $L$ with more than one element has a proper
congruence-preserving extension $K$ with the property that $L$
is an ideal in $K$.
 \end{theorem}

 \begin{proof}
 Take an element $a \in L$ such that $[a)$ (the dual ideal
generated by
$a$) has more than one element.  Then by
Lemma~\ref{L:basic}.(v), $A =
\new{[a)}$ is a proper congruence-preserving extension of
$[a)$ and $I = [a)$ is an ideal in $A$.  Now form the lattice
$K$ by gluing $L$ with the dual ideal $[a)$ to $A$ with the
ideal $I$.  It is clear that $K$ is a proper
congruence-preserving extension of $L$.
 \end{proof}

 \subsection*{Modularity and semimodularity}
 R. W. Quackenbush \cite{rQ85} proved that if $L$ is a modular
lattice, then  $M_3[L]$ is a semimodular lattice.  For our
construction, the analogous result fails: $\new{P}$ is not
semimodular, where $P$ is a projective plane (a modular
lattice).  Indeed, let $a$, $b$, $c$ be a triangle in $P$, with
sides $l$, $m$, $n$, that is, let $l$, $m$, $n$ be three
distinct lines in the plane $P$, and define the points $a = n
\mm m$, $b = n \mm l$, $c = m \mm l$.  Let $p$ be a point in
$P$ not on any one of these lines.  Then $\vv<p, \es, \es>$ is
an atom in $\new{P}$,  $\vv<a, b, c> \in \new{P}$ but
 \[
   \vv<\set{p}, \es, \es> \jj \vv<a, b, c> = \ol{\vv<p \jj a,
b, c>} =
      \vv<P, l, l>
 \] and
 \[
   \vv<a, b, c> < \vv<n, b, l> < \vv<P, l, l>,
 \] showing that $\new{P}$ is not semimodular.

Now we characterize when $\new{L}$ is modular.

\begin{theorem}\label{T:distr}
 Let $L$ be a lattice with more than one element.  Then
$\new{L}$ is modular iff $L$ is distributive.
 \end{theorem}
\begin{proof} If $L$ is distributive, then $\new{L} = M_3[L]$,
so  $\new{L}$ is modular by Theorem~\ref{T:Schmidt}.

Conversely, if $\new{L}$ is modular, then $L$ is modular since
it is a sublattice of $\new{L}$.  Now if $L$ is not
distributive, then $L$ contains an $M_3 = \set{o, a, b, c, i}$
as a sublattice.  By Lemma~\ref{L:basic}.(vi), the elements
 \[
    \vv<o, o, a>,
    \vv<o, c, a>,
    \vv<c, c, i>,
    \vv<i, i, i>,
    \vv<b, o, a>
 \] belong to $\new{L}$. Obviously,
 \[
   \vv<o, o, a> < \vv<o, c, a> < \vv<c, c, i> < \vv<i, i, i>
 \] and
 \[
   \vv<o, o, a> < \vv<b, o, a> < \vv<i, i, i>.
 \] To prove that these five elements form an $N_5$, it is
enough to prove that
 \[
   \vv<c, c, i> \mm \vv<b, o, a> = \vv<o, o, a>
 \] and
 \[
   \vv<o, c, a> \jj \vv<b, o, a> =  \vv<i, i, i>.
 \] The meet is obvious.  Now the join:
 \[
   \vv<o, c, a> \jj \vv<b, o, a> = \ol {\vv<b, c, a>} = \vv<i,
i, i>.
 \] So $\new{L}$ contains $N_5$ as a sublattice, contradicting
the assumption that $\new{L}$ is modular.  Therefore, $L$ is
distributive.
 \end{proof}

 \subsection*{Further results}
 $M_3[L]$ is not a lattice for a general $L$.  See, however, G.
Gr\"atzer and F. Wehrung \cite{GWa}, where a new concept of
$n$-modularity is introduced, for any natural number $n$.
Modularity is the same as $1$-modularity.

By definition, $n$-modularity is an identity; for larger $n$,
a weaker identity.  For an $n$-modular lattice $L$, $M_3[L]$ is
a lattice, a congruence-preserving extension of~$L$.

For distributive lattices (in fact, for $n$-modular lattices),
the construction $M_3[L]$ is a special case of the tensor
product construction of two semilattices with zero, see, for
instance, G. Gr\"atzer, H. Lakser, and R. W. Quackenbush
\cite{GLQ81} and R. W. Quackenbush \cite{rQ85}.  The $\new{L}$
construction is generalized in G. Gr\"atzer and F. Wehrung
\cite{GWb} to two bounded lattices; the new construction is
called \emph{box product}.  Some of the arguments of this
paper carry over to box products.

\section*{Problems}
\subsection*{Lattices}
As usual, let us denote by $\mbf{T}$, $\mbf{D}$, $\mbf{M}$, and
$\mbf{L}$ the variety of one-element, distributive, modular, and
all lattices, respectively. A variety $\mbf{V}$ is
\emph{nontrivial} if  $\mbf{V} \neq \mbf{T}$.

 Let us say that a variety $\mbf{V}$ of lattices has the
\emph{Congruence Preserving Extension Property} (CPEP, for short),
if every lattice in $\mbf{V}$ with more than one element has a
proper congruence-preserving extension in $\mathbf{V}$. It is
easy to see that no finitely generated lattice variety has CPEP.
(Indeed, by J\'onsson's lemma, a nontrivial finitely generated
lattice variety $\mbf{V}$ has a finite maximal subdirectly
irreducible member~$L$; if $K$ is a proper congruence-preserving
extension of $L$, then $K$ is also subdirectly irreducible and
$|L| >|K|$, a contradiction.)  In particular, $\mbf{D}$ does not
have CPEP.

Theorem~\ref{T:new} can be restated as follows: $\mbf{L}$ has
CPEP.

\begin{problem}
Find all lattice varieties $\mbf{V}$ with CPEP.  In particular,
does $\mbf{M}$ have CPEP?
\end{problem}

\subsection*{Groups}
  Let us say that a variety $\mbf{V}$ of groups has the
\emph{Normal Subgroup Preserving Extension Property} (NSPEP, for
short), if every group $G$ in $\mbf{V}$ with more than one
element has a proper supergroup $\ol{G}$ in $\mathbf{V}$ with the
following property: every normal subgroup $N$ in $G$ can be
uniquely represented in the form $\ol{N} \ii G$, where $\ol{N}$
is a normal subgroup of $\ol{G}$.

Not every group variety $\mbf{V}$ has NSPEP, for instance, the
variety $\mbf{A}$ of Abelian groups does not have NSPEP.

 \begin{problem}
Does the variety $\mbf{G}$ of all groups have NSPEP?  Find all
group varieties having NSPEP?
 \end{problem}

\subsection*{Rings}
 For ring varieties, we can similarly introduce the \emph{Ideal
Preserving Extension Property} (IPEP, for short).  The variety
$\mbf{R}$ of all (not necessarily commutative) rings has IPEP.
Indeed, if $R$ is a ring with more than one element, then embed
$R$ into $M_2(R)$ (the ring of $2 \times 2$ matrices over $R$)
with the diagonal map. The two-sided ideals of $M_2(R)$ are of the
form $M_2(I)$, where $I$ is a two-sided ideal of $R$, and $I =
M_2(I) \ii R$.

 \begin{problem}
Find all ring varieties having IPEP? In particular, does the
variety of all commutative rings have IPEP?
 \end{problem}

The second author found a positive answer for Dedekind domains:
every Dedekind domain with more than one element has a proper
ideal-preserving extension that is, in addition, a principal
ideal domain.

\section*{Acknowledgment} This work was partially completed
while the second author was visiting the University of
Manitoba. The excellent conditions provided by the Mathematics
Department, and, in particular, a quite lively seminar, were
greatly appreciated.

\end{document}